\documentclass{amsart}

\usepackage{amssymb}
\usepackage{mathtools}
\usepackage{mleftright}
\usepackage{smfthm}
\usepackage{stmaryrd}
\usepackage{tikz-cd}

\binoppenalty\maxdimen
\relpenalty\maxdimen

\title{Characterising properties of commutative rings using Witt vectors}
\author{Rubén Muñoz-\relax-Bertrand}
\address{Université Marie et Louis Pasteur, CNRS, LmB (UMR 6623), 25000 Besançon, France.}

\begin{document}

\begin{abstract}
We give equivalences between given properties of a commutative ring, and other properties on its ring of Witt vectors.
Amongst them, we characterise all commutative rings whose rings of Witt vectors are Noetherian.
We define a new category of commutative rings called preduced rings, and explain how it is the category of rings whose ring of Witt vectors has no $p$-torsion.
We then extend this characterisation to the torsion of the de Rham--Witt complex.
\end{abstract}

\maketitle

\section*{Introduction}

In this article $p$ a denotes prime number.
For any commutative ring $R$, we will denote by $W\mleft(R\mright)$ its ring of Witt vectors, by which we always mean $p$-typical Witt vectors.

We will denote by $\mleft[r\mright]\in W\mleft(R\mright)$ the multiplicative representative of $r\in R$, and respectively by $F$ and $V$ the Frobenius and Verschiebung maps on $W\mleft(R\mright)$.

\bigskip

Witt vectors play an important role in $p$-adic cohomology or in $p$-adic Hodge theory.
They have also found some applications to the theory of error correcting codes and cryptography.
It is thus natural to enquire about the properties of the ring $W\mleft(R\mright)$.

We shall be interested on how a given property of the commutative ring $R$ is equivalent to $W\mleft(R\mright)$ satisfying another property.
Many implications and equivalences were already known, sometimes only folklorically.
We recall them here and give the proofs when no reference has been found.

The main motivation for this study was to understand when a ring of Witt vectors was either Noetherian, or $p$-torsion free.
These two questions have arisen in the study of the de Rham--Witt complex.
We answer to them here.
We give a proof of the fact that Noetherian perfect rings of characteristic $p$ are the only commutative rings whose ring of Witt vectors is Noetherian.

For Witt vectors rings without $p$-torsion, we describe a new category of commutative rings called preduced rings.
This enables us to prove that a smooth commutative algebra over any commutative $\mathbb Z_{\mleft\langle p\mright\rangle}$-algebra is preduced if and only if its de Rham--Witt complex has no $p$-torsion.

To finish, we briefly study the open problem of coherent rings of Witt vectors.
We explain how, when $p$ is invertible in $R$, we retrieve the theory of uniformly coherent rings.

\bigskip

For the reader's comfort, the following tabular summarises all the equivalences which appear in this article.
On each line, a property of the commutative ring $R$ in the first column is equivalent to the property of $W\mleft(R\mright)$ in the second column.
\begin{center}
\begin{tabular}{|c|c|}\hline
\emph{$R$ is \ldots}&\emph{$W\mleft(R\mright)$ is \ldots}\\\hline
of characteristic $l$ coprime with $p$&of characteristic $l$ coprime with $p$\\\hline
of characteristic $pn$ for $n\in\mathbb N$&of characteristic $0$\\\hline
reduced&reduced\\\hline
$p$-adically complete and separated&$p$-adically complete and separated\\\hline
an integral domain of characteristic $p$&an integral domain\\\hline
a perfect cnd of characteristic $p$&a completely normal domain (cnd)\\\hline
Noetherian perfect of characteristic $p$&Noetherian\\\hline
a perfect field of characteristic $p$&a discrete valuation ring\\\hline
preduced (definition \ref{preduceddefinition})&$p$-torsion free\\\hline
\end{tabular}
\end{center}

\section*{Acknowledgments}

The author is under ``Contrat EDGAR-CNRS no 277952 UMR 6623, financé par la région Bourgogne-Franche-Comté.''.
This work has been financed by ANR-21-CE39-0009-BARRACUDA.
It was also partially written during my position at Université Paris-Saclay, UVSQ, CNRS, Laboratoire de mathématiques de Versailles, 78000, Versailles, France.

\section{Overture}

We assume the reader is acquainted with Witt vectors, and we refer to \cite{bourbakicommutativehuit} for an introduction.

We start with the well-known characterisation of the characteristic of the ring $W\mleft(R\mright)$.
As the author did not find any proof in the literature, it is given below for future reference.

After that, we shall study reduced and $p$-adically complete and separated rings.

\begin{enonce}{Proposition}\label{characteristiclcharacterisation}
Let $l\in\mathbb N^*$ be an integer coprime with $p$.
A commutative ring has characteristic $l$ if and only if its ring of Witt vectors has characteristic $l$.
\end{enonce}

\begin{proof}
The proposition is immediate when $l=1$, that is for the zero ring.

Otherwise, if $R$ has characteristic $l$, then $p$ is invertible in $R$.
Thus, by \cite[IX. \textsection1 proposition 2 and théorème 4]{bourbakicommutativehuit} we have an isomorphism of rings $W\mleft(R\mright)\cong R^\mathbb N$.

Conversely, if $W\mleft(R\mright)$ has characteristic $l$, then since $R\cong W\mleft(R\mright)/V\mleft(W\mleft(R\mright)\mright)$ as rings we see that the characteristic of $R$ must divide $l$, in particular it is coprime to $p$, so by our first point it must be $l$ itself.
\end{proof}

\begin{enonce}{Proposition}
A commutative ring has characteristic a multiple of $p$ if and only if its ring of Witt vectors has characteristic $0$.
\end{enonce}

\begin{proof}
Let $R$ be a commutative ring of characteristic $pn$, where $n\in\mathbb N$.
There is a surjective morphism of rings $W\mleft(R\mright)\to W\mleft(R/pR\mright)$, so for the implication it suffices to prove the case $n=1$.

The surjective morphism of rings $W\mleft(R\mright)\to R$ shows that the ring of Witt vectors of $R$ has characteristic a multiple of $p$.
But it is known that $p^u=V^u\mleft(1\mright)$ in $W\mleft(R\mright)$ for each $u\in\mathbb N$ \cite[IX. \textsection1 proposition 5]{bourbakicommutativehuit}, so we are done.

The converse follows from proposition \ref{characteristiclcharacterisation}.
\end{proof}

\begin{enonce}{Proposition}
A commutative ring is reduced if and only if its ring of Witt vectors also is.
\end{enonce}

\begin{proof}
Assume that $R$ is a non zero reduced commutative ring.
Then, $R$ has either prime or zero characteristic.
If its characteristic is not $p$, then $W\mleft(R\mright)$ is isomorphic to a subring of $R^\mathbb N$ \cite[IX. \textsection1 proposition 2 and théorème 4]{bourbakicommutativehuit}, so it is reduced.

If $R$ has characteristic $p$, then for any Witt vector $w\in W\mleft(R\mright)\smallsetminus V\mleft(W\mleft(R\mright)\mright)$ and any $n,u\in\mathbb N$ we have $V^u\mleft(w\mright)^n=V^{nu}\mleft(F^{\mleft(n-1\mright)u}\mleft(w^n\mright)\mright)$ \cite[IX. \textsection1 Proposition 5]{bourbakicommutativehuit}.
By definition, $w_0\neq 0$, therefore, $\mleft(w^n\mright)_0={w_0}^n\neq0$ by hypothesis on $R$.
But $F$ is injective, still by our hypothesis on $R$, and $V$ is always injective, so $W\mleft(R\mright)$ is reduced.

Conversely, assume that $W\mleft(R\mright)$ is reduced.
Then, for any $w\in R^*$ and any $n\in\mathbb N^*$ we have $\mleft[w\mright]^n\neq0$ so we conclude.
\end{proof}

\begin{enonce}{Proposition}\label{padicallycompleteandseparatedcharacterisation}
A commutative ring is $p$-adically complete and separated if and only if its ring of Witt vectors also is.
\end{enonce}

\begin{proof}
The implication is \cite[1. proposition 3]{thedisplayofaformal}.
Conversely, let $R$ be a commutative ring whose ring of Witt vectors is $p$-adically complete and separated.

Let $n\in\mathbb N^*$.
There is a unique morphism of rings $F\colon\mathbb Z\mleft[X_0,X_n\mright]\to\mathbb Z\mleft[X_0,X_n\mright]$ such that $F\mleft(X_0\mright)=\mleft(X_0-p^nX_n\mright)^p+p^n{X_n}^p$ and $F\mleft(X_n\mright)={X_n}^p$.
It is a lift of the Frobenius endomorphism of $\mathbb F_p\mleft[X_0,X_n\mright]$, promoting the $p$-torsion free ring $\mathbb Z\mleft[X_0,X_n\mright]$ to a $\delta$-ring \cite[remark 2.2]{prismsandprismaticcohomology}.

Furthermore, one notices that for the associated $\delta$-structure we get:

\begin{equation*}
\delta\mleft(X_0-p^nX_n\mright)=\frac{F\mleft(X_0\mright)}p-p^{n-1}F\mleft(X_n\mright)-\frac{\mleft(X_0-p^nX_n\mright)^p}p=0\text.
\end{equation*}

By \cite[lemma 2.9]{prismsandprismaticcohomology}, this implies that there exists a unique $\delta$-structure on $S_n\coloneqq\mathbb Z\mleft[X_0,X_n\mright]/\mleft\langle X_0-p^nX_n\mright\rangle$ such that the quotient map $\pi_n\colon\mathbb Z\mleft[X_0,X_n\mright]\to S_n$ is a morphism $\delta$-rings.

Let $r\in\bigcap_{n\in\mathbb N}p^nR$.
This means that for each $n\in\mathbb N^*$, there is $r_n\in R$ such that $p^nr_n=r$.
The evaluation map $e_n\colon\mathbb Z\mleft[X_0,X_n\mright]\to R$ satisfying $e_n\mleft(X_0\mright)=r$ and $e_n\mleft(X_n\mright)=r_n$ thus factors through $\pi_n$.
In turn, we find a morphism of $\delta$-rings $w_n\colon S_n\to W\mleft(R\mright)$ such that $\mleft(w_n\circ\pi_n\mleft(X_0\mright)\mright)_0=r$ and $\mleft(w_n\circ\pi_n\mleft(X_n\mright)\mright)_0=r_n$ \cite[théorème 4]{joyal}.

By construction, $w_n\circ\pi_n\mleft(X_0\mright)=p^nw_n\circ\pi_n\mleft(X_n\mright)$.
Thus $\lim_{n\to+\infty}w_n\circ\pi_n\mleft(X_0\mright)=0$ because we assumed that $W\mleft(R\mright)$ is $p$-adically complete and separated.
But the canonical map $W\mleft(R\mright)\to V\mleft(W\mleft(R\mright)\mright)$ is a continuous map for the $p$-adic topologies \cite[07E9]{stacksproject}, hence $r=0$.
\end{proof}

\section{Characterising domains}

In this section, we study how properties of integral domains behave with regard to the Witt vector functor.

\begin{enonce}{Lemma}\label{voneandp}
Let $R$ be a commutative ring.
Then, for any $w\in V\mleft(W\mleft(R\mright)\mright)$ we have:
\begin{equation*}
V\mleft(1\mright)w=pw\text.
\end{equation*}
\end{enonce}

\begin{proof}
We always have $V\mleft(1\mright)w=V\mleft(F\mleft(w\mright)\mright)$.
Let $x\in W\mleft(R\mright)$ satisfying $w=V\mleft(x\mright)$.
Thus, $V\mleft(1\mright)w=V\mleft(F\mleft(V\mleft(x\mright)\mright)\mright)$.
But $F\circ V$ is the multiplication by $p$, therefore $V\mleft(1\mright)w=pV\mleft(x\mright)=pw$.
\end{proof}

The previous lemma is only interesting when $p\neq0$ in $R$, otherwise $V\mleft(1\mright)=p$ and the lemma is vacuous.

\begin{enonce}{Proposition}\label{integraldomaincharacterisation}
A commutative ring is an integral domain of characteristic $p$ if and only if its ring of Witt vectors is an integral domain.
\end{enonce}

\begin{proof}
Assume that $R$ is an integral domain of characteristic $p$.
Then for any $w,x\in W\mleft(R\mright)\smallsetminus V\mleft(W\mleft(R\mright)\mright)$ and any $u,v\in\mathbb N$, applying \cite[IX. \textsection1 Proposition 5]{bourbakicommutativehuit} we get $V^u\mleft(w\mright)V^v\mleft(x\mright)=V^{u+v}\mleft(F^v\mleft(w\mright)F^u\mleft(x\mright)\mright)$.
Therefore, the $u+v$-th coefficient of that Witt vector is $w_0^{p^v}x_0^{p^u}$, which is not zero by hypothesis on $R$.
So $W\mleft(R\mright)$ is an integral domain.

Conversely, if $W\mleft(R\mright)$ is an integral domain, since it follows from lemma \ref{voneandp} that $V\mleft(1\mright)\mleft(V\mleft(1\mright)-p\mright)=0$, we must have $V\mleft(1\mright)=p$ which only happens when $R$ has characteristic $p$.
Then, $R$ is an integral domain because $\mleft[\bullet\mright]$ is an injective and multiplicative map.
\end{proof}

We now turn to normal domains.
These objects are sometimes known as integrally closed domains, but here we follow the terminology of \cite[0309]{stacksproject}.

For once, we only give one implication.
The characterisation of all commutative rings whose ring of Witt vectors is a normal domain is a more difficult question, which was the motivation of \cite{onthewittvectorsofperfectringsinpositivecharacteristic}, and for which only partial results are known.

\begin{enonce}{Lemma}\label{normalimplication}
If a commutative ring has a ring of Witt vectors which is a normal domain, then it is a perfect integral domain of characteristic $p$.
\end{enonce}

\begin{proof}
Assume that $R$ is a commutative ring such that $W\mleft(R\mright)$ is a normal domain.
By proposition \ref{integraldomaincharacterisation}, $R$ is an integral domain of characteristic $p$, in particular it is reduced.

Let $r\in R$.
Consider the polynomial $X^p-\mleft[r\mright]$.
In $\operatorname{Frac}\mleft(W\mleft(R\mright)\mright)$, applying \cite[IX. \textsection1 Proposition 5]{bourbakicommutativehuit} we see that it has a root because: 
\begin{equation*}
\mleft(\frac{V\mleft(\mleft[r\mright]\mright)}p\mright)^p=\frac{p^{p-1}V\mleft(\mleft[r\mright]^p\mright)}{p^p}=\mleft[r\mright]\text.
\end{equation*}

We assumed that $W\mleft(R\mright)$ is a normal domain, so this means that $p$ divides $V\mleft(\mleft[r\mright]\mright)$ in $W\mleft(R\mright)$.
But this is only possible if $\mleft[r\mright]$ is in the image of $F$, so $R$ is semiperfect, hence perfect.
\end{proof}

The last characterisation of this section is about completely normal domains.
They are sometimes called completely integrally closed domains, again we rather follow the terminology in \cite[00GW]{stacksproject}.

\begin{enonce}{Proposition}
A commutative ring is a perfect integral domain of characteristic $p$ if and only if its ring of Witt vectors is a completely normal domain.
\end{enonce}

\begin{proof}
The implication is \cite[corollary 5.10]{onthewittvectorsofperfectringsinpositivecharacteristic}, so let us focus on the converse.

Assume that $R$ is a commutative ring such that $W\mleft(R\mright)$ is a completely normal domain.
By \cite[00GX]{stacksproject} and lemma \ref{normalimplication}, $R$ is a perfect integral domain of characteristic $p$.
Let $x\in\operatorname{Frac}\mleft(R\mright)$ such that there exists $r\in R$ such that for each $n\in\mathbb N$ we have $rx^n\in R$.

Consider the commutative diagram of commutative rings, in which every arrow is injective:
\begin{equation*}
\begin{tikzcd}
W\mleft(R\mright)\arrow[r]\arrow[d]&W\mleft(\operatorname{Frac}\mleft(R\mright)\mright)\arrow[d]\\
\operatorname{Frac}\mleft(W\mleft(R\mright)\mright)\arrow[r,"\psi"]&\operatorname{Frac}\mleft(W\mleft(\operatorname{Frac}\mleft(R\mright)\mright)\mright)\text.
\end{tikzcd}
\end{equation*}

The image of $\mleft[x\mright]$ in $\operatorname{Frac}\mleft(W\mleft(\operatorname{Frac}\mleft(R\mright)\mright)\mright)$ is also in the image of $\psi$ because $\mleft[\bullet\mright]$ is multiplicative.
For the same reason, we see that $\mleft[x\mright]\mleft[r\mright]^n$ is the image of an element in $W\mleft(R\mright)$ for every $n\in\mathbb N$, so $\mleft[x\mright]$ comes from $W\mleft(R\mright)$, that is, $R$ is a completely normal domain.
\end{proof}

\section{Characterisation of Noetherian Witt vectors rings}

The main point of this section is to show that there are commutative rings having a Noetherian ring of Witt vectors which are not discrete valuation rings.

\begin{enonce}{Lemma}\label{pmcharwittnoether}
A commutative ring whose ring of Witt vectors is Noetherian has characteristic $p^{m}$ for some $m\in\mathbb{N}$.
\end{enonce}

\begin{proof}
Let $R$ be a commutative ring with characteristic $c\in\mathbb{N}\smallsetminus\mleft\{p^{m}\mid m\in\mathbb{N}\mright\}$.
Consider $I\coloneqq\mleft\{r\in R\mid\exists m\in\mathbb{N},\ p^{m}r=0\mright\}$. By hypothesis on $R$, it is a proper ideal of $R$.
By construction $R/I$ is a $p$-torsion-free ring.
In particular, if it has not characteristic $0$, then $p$ is a unit in $R/I$.

Assume that $W\mleft(R\mright)$ is Noetherian.
We deduce from the surjective ring morphism $W\mleft(R\mright)\to W\mleft(R/I\mright)$ that $W\mleft(R/I\mright)$ is Noetherian \cite[00FN]{stacksproject}.
If $p$ is a unit in $R/I$, then we have an isomorphism of rings $W\mleft(R/I\mright)\cong\mleft(R/I\mright)^{\mathbb{N}}$ \cite[IX. \textsection1 proposition 2]{bourbakicommutativehuit}.
In particular, $W\mleft(R/I\mright)$ cannot be Noetherian so we get to a contradiction, so $R$ must have characteristic $0$.

But if one takes $l$ a prime number distinct from $p$, then $W\mleft(R/lR\mright)$ must be Noetherian, which is impossible by the above argument.
\end{proof}

\begin{enonce}{Lemma}\label{prwittnoether}
Any commutative ring $R$ such that $W\mleft(R\mright)$ is Noetherian satisfies $pR=\mleft\{0\mright\}$.
\end{enonce}

\begin{proof}
We already know that there exists $m\in\mathbb{N}$ such that $R$ has characteristic $p^{m}$, because of lemma \ref{pmcharwittnoether}.
We only have to show that $m$ is not greater than $2$. If it is the case, then $p^{m-1}$ is not zero, but $\mleft(p^{m-1}\mright)^{2}=0$.

Therefore, using \cite[IX. \textsection1 proposition 4]{bourbakicommutativehuit} we get the following set inclusion for all $u\in\mathbb{N}^*$:
\begin{equation*}
\mleft\langle V^{u}\mleft(\mleft[p^{m-1}\mright]\mright)\mright\rangle\subset\mleft\{V^{u}\mleft(\mleft[r\mright]\mright)\mid r\in R\mright\}\text{,}
\end{equation*}
where the left-hand side is the ideal of $W\mleft(R\mright)$ generated by $V^{u}\mleft(\mleft[p^{m-1}\mright]\mright)$. Applying the cited proposition again, we deduce that the ideal of $W\mleft(R\mright)$ generated by the set $\mleft\{V^{u}\mleft(\mleft[p^{m-1}\mright]\mright)\mid u\in\mathbb{N}^{*}\mright\}$ cannot be finitely generated.
\end{proof}

\begin{enonce}{Proposition}\label{noetherianwitt}
Let $R$ be a commutative ring. Then $W\mleft(R\mright)$ is Noetherian if and only if $R$ is either the zero ring, or a Noetherian perfect ring of characteristic $p$.
\end{enonce}

\begin{proof}
We will first assume that $W\mleft(R\mright)$ is Noetherian. After having applied lemma \ref{prwittnoether} to $R$, we only have to show that if $R$ has characteristic $p$, then it is Noetherian and perfect.
Of course the quotient ring $R\cong W\mleft(R\mright)/V\mleft(W\mleft(R\mright)\mright)$ must be Noetherian.
Moreover, a Noetherian ring with characteristic $p$ is perfect if and only if it is semiperfect \cite[06RN]{stacksproject}, so we only have to show that $R$ is semiperfect.

If it is not, then we can find $r\in R$ which is not in the image of the Frobenius endomorphism. Consider $u\in\mathbb{N}^{*}$ and $w\in W\mleft(R\mright)$.
We can write $w=\mleft[w_{0}\mright]+V\mleft(x\mright)$ for some $x\in W\mleft(R\mright)$, so that $F^{u}\mleft(w\mright)=\mleft[{w_{0}}^{p^{u}}\mright]+pF^{u-1}\mleft(x\mright)$.
Hence:
\begin{equation*}
V^{u}\mleft(\mleft[r\mright]\mright)w=V^u\mleft(\mleft[r{w_0}^{p^u}\mright]\mright)+pV^u\mleft(F^{u-1}\mleft(x\mright)\mright)\text.
\end{equation*}

We can deduce from this that $V^{u+1}\mleft(\mleft[r\mright]\mright)$ does not belong to the ideal of $W\mleft(R\mright)$ generated by the set $\mleft\{V^{m}\mleft(\mleft[r\mright]\mright)\mid m\in\mleft\llbracket1,u\mright\rrbracket\mright\}$.
Indeed, if for all $m\in\mleft\llbracket1,u\mright\rrbracket$ there exists a $w\mleft(m\mright)\in W\mleft(R\mright)$ such that $V^{u+1}\mleft(\mleft[r\mright]\mright)=\sum_{m=1}^{u}V^{m}\mleft(\mleft[r\mright]\mright)w\mleft(m\mright)$, then we get after reduction modulo $p$:
\begin{equation*}
\overline{V^{u+1}\mleft(\mleft[r\mright]\mright)}=\overline{\sum_{m=1}^{u}V^{m}\mleft(\mleft[r{{w\mleft(m\mright)}_{0}}^{p^{u}}\mright]\mright)}\text{.}
\end{equation*}

But this is impossible, because it implies that $r$ is a power of $p$ \cite[IX. \textsection1 proposition 5]{bourbakicommutativehuit}.
So $R$ must be semiperfect, hence perfect.

Conversely, if $R$ is a Notherian perfect ring of characteristic $p$, then $W\mleft(R\mright)$ is $p$-adically complete and separated, and $W\mleft(R\mright)/pW\mleft(R\mright)\cong R$ as rings \cite[IX. \textsection1 proposition 7]{bourbakicommutativehuit}.
This implies that $W\mleft(R\mright)$ is Noetherian \cite[05GH]{stacksproject}.
\end{proof}

The following proposition is probably well-known to experts, but the author could not find a proof in the literature.

\begin{enonce}{Proposition}
A commutative ring is a perfect field with characteristic $p$ if and only if its ring of Witt vectors is a discrete valuation ring.
\end{enonce}

\begin{proof}
Let $R$ be a commutative ring such that $W\mleft(R\mright)$ is a discrete valuation ring.
By definition, $W\mleft(R\mright)$ is an integral domain, and it is Noetherian \cite[00II]{stacksproject}.
We have shown in proposition \ref{integraldomaincharacterisation} that $R$ must be an integral domain of characteristic $p$, and in proposition \ref{noetherianwitt} that it also has to be perfect.

By proposition \ref{padicallycompleteandseparatedcharacterisation}, $W\mleft(R\mright)$ is $p$-adically complete and separated.
It thus follows from \cite[III. \textsection2 lemme 3]{algebrecommutativeun} that $p$ is in the Jacobson radical of $W\mleft(R\mright)$.

By \cite[00PD]{stacksproject}, $W\mleft(R\mright)$ has a unique maximal ideal generated by some $w\in{W\mleft(R\mright)}^*$ such that we can find $u\in{W\mleft(R\mright)}^\times$ and $n\in\mathbb N$ such that $p=w^nu$.
Since $R$ has characteristic $p$, this implies that ${w_0}^nu_0=0$.
Since $u_0$ must be a unit, as $R$ is an integral domain this implies that $w_0=0$.

It follows from the perfection of $R$ and from \cite[IX. \textsection1 proposition 7]{bourbakicommutativehuit} that $p$ divides $w$; in other words, $pW\mleft(R\mright)$ is the unique maximal ideal of $W\mleft(R\mright)$.
Thus, $R$ is a field.

The converse is well documented \cite[IX. \textsection1 proposition 8]{bourbakicommutativehuit}.
\end{proof}

The last two propositions imply that there exist Noetherian Witt vectors rings which are not discrete valuation rings. Such an example is the ring of Witt vectors over $\mathbb{F}_{p}\mleft[X\mright]/\mleft\langle X^{p}-X\mright\rangle$, which is a Noetherian perfect ring of positive characteristic, but not a field.

\section{Preduced rings}

We now turn to the characterisation of all commutative rings $R$ whose ring of Witt vectors has no $p$-torsion.
At the end of this section, we will explain how this characterisation can be extended to the de Rham--Witt complex.

It turns out that the corresponding property on $R$ is not a usual one, so we will introduce a new category of commutative preduced rings.
To define it, we will employ the $V$-adic pseudovaluation of the ring of Witt vectors of a commutative ring $R$:
\begin{equation*}
\operatorname v_V\colon\begin{array}{rl}W\mleft(R\mright)\to&\mathbb{R}\cup\mleft\{+\infty,-\infty\mright\}\\w\mapsto&\max\mleft\{m\in\mathbb{N}^{*}\mid w\in V^{m}\mleft(W\mleft(R\mright)\mright)\mright\}\end{array}\text.
\end{equation*}

\begin{enonce}{Proposition}\label{wittptorsion}
On a commutative ring $R$, the following three conditions are equivalent:
\begin{align*}
\forall w\in{W\mleft(R\mright)}^{*},\ &pw\neq0\text,\\
\forall w\in W\mleft(R\mright),\ &\operatorname{v}_{V}\mleft(w\mright)\leqslant\operatorname{v}_{V}\mleft(pw\mright)\leqslant\operatorname{v}_{V}\mleft(w\mright)+1\text,\\
\forall r\in R^{*},\ &\mleft(pr=0\mright)\implies\mleft(r^{p}\neq0\mright)\text.
\end{align*}
\end{enonce}

\begin{proof}
Assume that $W\mleft(R\mright)$ has no $p$-torsion. Let $r\in R^{*}$ be such that $pr=0$. Since $\mleft[r\mright]\neq0$ we also have $p\mleft[r\mright]=F\mleft(V\mleft(\mleft[r\mright]\mright)\mright)\neq0$ by hypothesis. As $r$ is $p$-torsion, we have $F\mleft(V\mleft(\mleft[r\mright]\mright)\mright)=\mleft(0,r^{p},0,0,\ldots\mright)$ according to \cite[IX. \textsection1 N$^\circ$ 3]{bourbakicommutativehuit}. So $r^{p}\neq0$.

Assume now that $R$ is a commutative ring such that all $r\in R^{*}$ with $pr=0$ satisfy $r^{p}\neq0$. Let $w\in W\mleft(R\mright)$. If $w\neq0$, then we can find $m\in\mathbb{N}$ and $y\in W\mleft(R\mright)\smallsetminus V\mleft(W\mleft(R\mright)\mright)$ such that $V^{m}\mleft(y\mright)=w$. We have $pw=V^{m}\mleft(py\mright)$. Using \cite[IX. \textsection1 N$^\circ$ 3]{bourbakicommutativehuit} again, we get that $\mleft(py\mright)_{0}=py_{0}$. If $py_{0}\neq0$, then $\operatorname{v}_{V}\mleft(py\mright)=\operatorname{v}_{V}\mleft(y\mright)$.

Otherwise, then $py_{0}=0$, which implies that ${y_{0}}^{p}\neq0$. Still according to Bourbaki we have $\mleft(py\mright)_{1}={y_{0}}^{p}+py_{1}$. To obtain $\operatorname{v}_{V}\mleft(py\mright)=\operatorname{v}_{V}\mleft(y\mright)+1$, we thus need to show that ${y_{0}}^{p}+py_{1}\neq0$. If it was not the case, we would get:
\begin{equation*}
{y_{0}}^{p^{2}}=-\mleft(-p\mright)^{p-1}py_{1}{y_{1}}^{p-1}=-\mleft(-p\mright)^{p-2}py_{0}\mleft(y_{0}y_{1}\mright)^{p-1}=0\text{.}
\end{equation*}

In particular, by hypothesis on $R$ we must have $p{y_{0}}^{p}\neq0$, which is impossible. So we have $\operatorname{v}_{V}\mleft(py\mright)=\operatorname{v}_{V}\mleft(y\mright)+1$. The case $w=0$ is obvious.

It is clear that the condition on the $V$-adic pseudovaluation implies that $W\mleft(R\mright)$ has no $p$-torsion.
\end{proof}

\begin{enonce}{Definition}\label{preduceddefinition}
We shall call a commutative ring $R$ satisfying any of the equivalent conditions in proposition \ref{wittptorsion} \emph{preduced} (short for ``$p$-reduced'' or ``pre-reduced'').
\end{enonce}

Obviously, all reduced commutative rings, a fortiori all integral perfectoid rings, are preduced.

In characteristic $p$, being preduced is equivalent to being reduced.

We shall denote by $\mathsf{PCRing}$ the full subcategory in the category $\mathsf{CRing}$ of commutative rings whose objects are preduced rings.

\begin{enonce}{Lemma}\label{preducedinduction}
A commutative ring $R$ is preduced if and only if for every $r\in R^*$ we have:
\begin{equation*}
\forall n\in\mathbb N^*,\ \mleft(p^nr\neq0\mright)\vee\mleft(r^{pn}\neq0\mright)\text.
\end{equation*}

This condition is equivalent to:
\begin{equation*}
\mleft(\forall n\in\mathbb N^*,\ p^nr\neq0\mright)\vee\mleft(\forall n\in\mathbb N^*,\ r^{pn}\neq0\mright)\text.
\end{equation*}
\end{enonce}

\begin{proof}
We prove the first equivalence by induction on $n\in\mathbb N^*$, the case $n=1$ being the definition of a preduced ring (thus giving the converse).

Assume that the lemma is shown for some $n\in\mathbb N^*$.
If $p^nr\neq0$, then either $p^{n+1}r\neq0$, or $p^{pn}r^p\neq0$.
But the second possibility implies the first one.

If $r^{pn}\neq0$, then either $pr^{pn}\neq0$ or $r^{p^2n}\neq0$.
The second implication would prove the lemma, so let us see what happens when it fails.
Then, we must have $pr^{pn}\neq0$, so either $p^2r^{pn}\neq0$ or $p^pr^{p^2n}\neq0$.
The latter can not happen, so with this reasoning we find by induction on $i\in\mathbb N$ that we always have $p^ir^{pn}\neq0$, in particular $p^{n+1}r\neq0$.

The above argument gives the implication between the two centred properties, whose converse is immediate.
\end{proof}

Notice in particular that a non-zero preduced commutative $\mathbb Z_{\langle p\rangle}$-algebra is therefore either of characteristic $0$, or a reduced ring of characteristic $p$.

\begin{enonce}{Proposition}\label{universalpreduced}
For any commutative ring $R$, let:
\begin{equation*}
\mathfrak I\mleft(R\mright)\coloneqq\mleft\{r\in R\mid\exists n\in\mathbb N^*,\ \mleft(p^nr=0\mright)\wedge\mleft(r^{pn}=0\mright)\mright\}\text.
\end{equation*}

It is an ideal of $R$, and $R/\mathfrak I\mleft(R\mright)$ is a preduced ring.

Moreover, for any morphism of commutative rings $R\to S$ there is a canonical morphism of preduced rings $R/\mathfrak I\mleft(R\mright)\to S/\mathfrak I\mleft(S\mright)$, naturally defining a left adjoint functor to the obvious forgetful functor $\operatorname{For}\colon\mathsf{PCRing}\to\mathsf{CRing}$.

More colloquially, $R\to R/\mathfrak{I}\mleft(R\mright)$ is universal for maps to preduced rings.
\end{enonce}

\begin{proof}
Let $R$ be a commutative ring.
Let $r,s\in\mathfrak I\mleft(R\mright)$.
Then we can find $\mathbb N^*$ large enough so that $r^{pn}=0$, $s^{pn}=0$, $p^nr=0$ and $p^ns=0$.
Of course, $p^n\mleft(r-s\mright)=0$.
For every integer $l$ between $0$ and $2n$ we have $r^{2pn-l}s^l=0$.
In particular $\mleft(r-s\mright)^{2pn}=0$, from which it is easy to see that $\mathfrak I\mleft(R\mright)$ is indeed an ideal of $R$.

Let $r\in R\smallsetminus\mathfrak I\mleft(R\mright)$.
If $pr\in\mathfrak I\mleft(R\mright)$, then $p^{n+1}r=0$ for all large enough $n\in\mathbb N^*$, so by hypothesis on $r$ we must have $r^{p\mleft(n+1\mright)}\neq0$.
This holds in particular when $p$ divides $n+1$, which implies that $r^p\notin\mathfrak I\mleft(R\mright)$, meaning that $R/\mathfrak I\mleft(R\mright)$ is preduced.

Let $P$ be a preduced commutative ring, and let $\varphi\colon R\to P$ be a morphism of rings.
Let $r\in\mathfrak I\mleft(R\mright)$, and let $n\in\mathbb{N}^*$ be an integer such that $r^{pn}=0$ and $p^nr=0$.
Thus, $\varphi\mleft(r\mright)^{pn}=0$ and $p^n\varphi\mleft(r\mright)=0$.
By lemma \ref{preducedinduction}, this means that $\varphi\mleft(r\mright)=0$.

Composition with the canonical morphism $\bullet\to\bullet/\mathfrak I\mleft(\bullet\mright)$ thus lead to an isomorphism of bifunctors:
\begin{equation*}
\operatorname{Hom}_{\mathsf{PCRing}}\mleft(\bullet/\mathfrak I\mleft(\bullet\mright),\bullet\mright)\cong\operatorname{Hom}_{\mathsf{CRing}}\mleft(\bullet,\operatorname{For}\mleft(\bullet\mright)\mright)
\end{equation*}
whose functoriality is straightforward, and left to the careful reader.
\end{proof}

\begin{enonce}{Proposition}\label{localisationpreduced}
Any localisation of a preduced commutative ring is preduced.
\end{enonce}

\begin{proof}
Let $S$ be a multiplicative subset of a commutative preduced ring $R$. Let $r\in R$ and $s\in S$. Assume that $\frac{r}{s}\neq0$ and $p\frac{r}{s}=0$ in $R\mleft[S^{-1}\mright]$. Then there exists $t\in S$ such that $tpr=0$. Assume also that $ur^{p}=0$ for some $u\in S$, then $\mleft(tur\mright)^{p}=0$ and $ptur\neq0$ as $R$ is preduced, but this is impossible. In particular, $\mleft(\frac{r}{s}\mright)^{p}\neq0$.
\end{proof}

\begin{enonce}{Proposition}\label{spectrumpreduced}
Let $R$ be a commutative ring.
Then the following conditions are equivalent:
\begin{enumerate}
\item the ring $R$ is preduced;
\item for all open sets $U$ in $\operatorname{Spec}\mleft(R\mright)$ the ring $\Gamma\mleft(U,\mathcal{O}_{\operatorname{Spec}\mleft(R\mright)}\mright)$ is preduced;
\item for all $\mathfrak{p}\in\operatorname{Spec}\mleft(R\mright)$ the stalk of $\mathcal{O}_{\operatorname{Spec}\mleft(R\mright)}$ at $\mathfrak{p}$ is reduced.
\end{enumerate}
\end{enonce}

\begin{proof}
This is standard hocus-pocus. Let $R$ be a commutative ring, such that the stalks on $\operatorname{Spec}\mleft(R\mright)$ are preduced. Let $U$ be an open set of $\operatorname{Spec}\mleft(R\mright)$. Let $x\in{\Gamma\mleft(U,\mathcal{O}_{\operatorname{Spec}\mleft(R\mright)}\mright)}^{*}$ satisfying $px=0$. Then there exists $\mathfrak{p}\in U$ such that the image $\widetilde{x}$ of $x$ in the stalk at $\mathfrak{p}$ is not zero, so $\widetilde{x}^{p}\neq0$ and $x^{p}\neq0$ by \cite[0079]{stacksproject}.

The converse is given by proposition \ref{localisationpreduced}.
\end{proof}

We will now study the behaviour of preduced rings with the de Rham--Witt complex.
Given a morphism $k\to R$ of commutative $\mathbb Z_{\mleft\langle p\mright\rangle}$-algebras, we shall denote by $W\Omega_{R/k}$ the associated de Rham--Witt complex.
For any $n\in\mathbb N^*$, its truncated de Rham--Witt complex will be denoted by $W_n\Omega_{R/k}$.
For the definitions, we refer to \cite{derhamwittcohomologyforaproperandsmoothmorphism}.

\begin{enonce}{Lemma}\label{etalepreduced}
Let $m\in\mathbb{N}^{*}$. Assume that $k$ is a preduced commutative $\mathbb{Z}_{\langle p\rangle}$-algebra. Let $R$ be a commutative ring such that there exists $n\in\mathbb{N}$ and an étale morphism of rings $k\mleft[X_{1},\ldots,X_{n}\mright]\to R$. Then the canonical projection $W_{m+1}\Omega_{R/k}\to W_{m}\Omega_{R/k}$ factors through a morphism of $W_{m+1}\mleft(k\mright)$-modules:
\begin{equation*}
W_{m+1}\Omega_{R/k}/\mleft\{w\in W_{m+1}\Omega_{R/k}\mid pw=0\mright\}\to W_{m}\Omega_{R/k}\text{.}
\end{equation*}
\end{enonce}

\begin{proof}
When $R=k\mleft[X_{1},\ldots,X_{n}\mright]$ for some $n\in\mathbb{N}$, the proposition is a consequence of \cite[corollary 2.13]{derhamwittcohomologyforaproperandsmoothmorphism}.
Indeed, if $w\in W_{m+1}\mleft(R\mright)^{*}$ satisfies $pw=0$, then by virtue of proposition \ref{wittptorsion} we must have $\operatorname{v}_{V}\mleft(w\mright)=m$.

On account of \cite[theorem 10.4.]{integralpadichodgetheory}, the canonical morphism of $W_{m}\mleft(k\mright)$-algebras $W_{m}\mleft(k\mleft[X_{1},\ldots,X_{n}\mright]\mright)\to W_{m}\mleft(R\mright)$ is étale.
So, by the étale base change property of the de Rham--Witt complex \cite[lemma 10.8.]{integralpadichodgetheory}, we have an isomorphism of $W_{m}\mleft(R\mright)$-modules:
\begin{multline*}
\mleft\{w\in W_{m+1}\Omega_{R/k}\mid pw=0\mright\}\\\cong\mleft\{w\in W_{m+1}\Omega_{k\mleft[X_{1},\ldots,X_{n}\mright]/k}\mid pw=0\mright\}\otimes_{W_{m}\mleft(k\mleft[X_{1},\ldots,X_{n}\mright]\mright)}W_{m}\mleft(R\mright)\text{,}
\end{multline*}
and the general case follows.
\end{proof}

\begin{enonce}{Proposition}\label{surjectivedrwiso}
Let $n\in\mathbb{N}^{*}$. Let $k$ be a commutative $\mathbb{Z}_{\mleft\langle p\mright\rangle}$-algebra. Let $I$ be an ideal of $k$. Let $R$ be a commutative $k$-algebra. Then we have an isomorphism of $W_{n}\mleft(k\mright)$-dgas:
\begin{equation*}
W_{n}\Omega_{\mleft(R/IR\mright)/k}\cong W_{n}\Omega_{\mleft(R/IR\mright)/\mleft(k/I\mright)}\text{.}
\end{equation*}

Moreover, the kernel of the natural morphism $W_{n}\Omega_{R/k}\to W_{n}\Omega_{\mleft(R/IR\mright)/\mleft(k/I\mright)}$ is the differential graded ideal generated by $\operatorname{Ker}\mleft(W_{n}\mleft(R\mright)\to W_{n}\mleft(R/IR\mright)\mright)$.
\end{enonce}

\begin{proof}

The first part of the statement follows from the same reasoning than in \cite[proposition 1.8]{derhamwittcohomologyforaproperandsmoothmorphism}, except that here our morphism is surjective instead of being unramified, so that the hypothesis on $R$ unneeded.

From \cite[lemma 10.9]{integralpadichodgetheory}, we deduce the second part of the statement.
\end{proof}

\begin{enonce}{Theorem}
Assume that $k$ is a commutative $\mathbb Z_{\mleft\langle p\mright\rangle}$-algebra.
Let $R$ be a smooth commutative $k$-algebra.
Then $W\Omega_{R/k}$ has no $p$-torsion if and only if $R$ is preduced.
\end{enonce}

\begin{proof}
The condition on $R$ is necessary by proposition \ref{wittptorsion}. Conversely, assume that $R$ is preduced. We can then assume that $k$ is also preduced by virtue of propositions \ref{surjectivedrwiso} and \ref{universalpreduced}.

Let $w\in W\Omega_{R/k}$ such that $pw=0$. Assume that there exists $m\in\mathbb{N}^{*}$ such that the image of $w$ in $W_{m}\Omega_{R/k}$ is not zero; that is, assume that $w\neq0$.

By \cite[039Q]{stacksproject}, we know we can find $n\in\mathbb{N}$, as well as two commutative rings $R'$ and $k'$ and a commutative diagram of affine schemes:
\begin{equation*}
\begin{tikzcd}
\operatorname{Spec}\mleft(R\mright)\arrow[rr]&&\operatorname{Spec}\mleft(k\mright)\\
\operatorname{Spec}\mleft(R'\mright)\arrow[u]\arrow[r,"\psi"]&\operatorname{Spec}\mleft(k'\mleft[X_{1},\ldots,X_{n}\mright]\mright)\arrow[r]&\operatorname{Spec}\mleft(k'\mright)\arrow[u]
\end{tikzcd}
\end{equation*}
where $\psi$ is étale, the vertical arrows are open immersions, and such that the image of $w$ through $W_{m}\Omega_{R/k}\to W_{m}\Omega_{R'/k'}$ is also not zero (for more details about the sheaf property of the truncated de Rham--Witt complex, see \cite[(1.36)]{derhamwittcohomologyforaproperandsmoothmorphism} and the following discussion, as well as their appendix).

Notice that $k'$ is preduced: this is proposition \ref{spectrumpreduced}.
Therefore we can apply lemma \ref{etalepreduced}, and as $pw=0$ we get that the image of $w$ in $W_{m}\Omega_{R'/k'}$ should be zero, which contradicts our hypothesis $w\neq0$.

\end{proof}

\section{Coherence}

A more profound question, generalising \cite[question 3]{onthewittvectorsofperfectringsinpositivecharacteristic}, is the following: which commutative rings have a coherent ring of Witt vectors?

There are a few results for this question, when $R$ is a valuation ring of characteristic $p$.
For instance, when the valuation of $R$ is non trivial and its value group is not isomorphic to $\mathbb R$, then $W\mleft(R\mright)$ is not coherent \cite[theorem 1.2]{someringtheoreticpropertiesof}.
This is neither the case if we assume instead that the valuation of $R$ has rank $1$ and is non discrete, and that $R$ is complete with respect to its valuation \cite[theorem 1.1]{langludwig}.

As we are about to see, the case of rings in which $p$ is invertible is much easier.

First, we recall that a uniformly coherent commutative ring is a commutative ring $R$ such that there exists a map $\phi\colon\mathbb N\to\mathbb N$ such that for each $n\in\mathbb N$, the kernel of every morphism of $R$-modules $R^n\to R$ can be generated by at most $\phi\mleft(n\mright)$ elements.

Of course, every uniformly coherent commutative ring is coherent.
However, not every coherent commutative ring is uniformly coherent.
In fact, there are even examples of Noetherian commutative rings which are not uniformly coherent \cite[proposition 3.1]{surluniformecoh}.

\begin{enonce}{Proposition}
Let $R$ be a commutative ring in which $p$ is invertible.
Then, $W\mleft(R\mright)$ is coherent if and only if $R$ is uniformly coherent.
\end{enonce}

\begin{proof}
By assumption on $R$, we have $W\mleft(R\mright)\cong R^\mathbb N$ as rings \cite[IX. \textsection1 proposition 2 and théorème 4]{bourbakicommutativehuit}.
But in that setting, the equivalence is already known \cite[théorème 16]{anneauxuniform}.
\end{proof}


\begin{thebibliography}{Bou06b}
\bibitem[BMS]{integralpadichodgetheory}Bhargav \textsc{Bhatt}, Matthew \textsc{Morrow}, Peter \textsc{Scholze}, \textit{Integral $p$-adic Hodge theory}. Publications mathématiques de l'IHÉS \textbf{128} (2018), pp. 219--397.

\bibitem[BS]{prismsandprismaticcohomology}Bhargav \textsc{Bhatt}, Peter \textsc{Scholze}, \textit{Prisms and prismatic cohomology}. Annals of Mathematics \textbf{196} n°3 (2022), pp. 1135--1275.

\bibitem[Bou06a]{algebrecommutativeun}Nicolas \textsc{Bourbaki}, \textit{Algèbre commutative : Chapitres 1 à 4}. Springer-Verlag (2006).

\bibitem[Bou06b]{bourbakicommutativehuit}Nicolas \textsc{Bourbaki}, \textit{Algèbre commutative : Chapitres 8 et 9}. Springer-Verlag (2006).

\bibitem[Joy]{joyal}André \textsc{Joyal}, \textit{$\delta$-anneaux et vecteurs de Witt}. Mathematical Reports of the Academy of Science of the Royal Society of Canada \textbf{7} n°3 (1985), pp. 177--182.

\bibitem[Ked]{someringtheoreticpropertiesof}Kiran Sridhara \textsc{Kedlaya}, \textit{Some Ring-Theoretic Properties of $\mathbf A_\mathrm{inf}$}. Simons Symposia \textbf{$p$-adic Hodge theory} (2020), pp. 129--141.

\bibitem[LL]{langludwig}Jaclyn \textsc{Lang}, Judith \textsc{Ludwig}, \textit{$\mathbb A_\mathrm{inf}$ is infinite dimensional}. Journal of the Institute of Mathematics of Jussieu \textbf{20} n°6 (2020), pp. 1983--1989.

\bibitem[LZ]{derhamwittcohomologyforaproperandsmoothmorphism}Andreas \textsc{Langer}, Thomas \textsc{Zink}, \textit{De Rham--Witt cohomology for a proper and smooth morphism}. Journal of the Institute of Mathematics of Jussieu \textbf{3} n°2 (2004), pp. 231--314.

\bibitem[Que]{surluniformecoh}Yann \textsc{Quentel}, \textit{Sur l'uniforme cohérence des anneaux noethériens}. Comptes rendus hebdomadaires des séances de l'Académie des Sciences Série A \textbf{275} n°16 (1972), pp. 753--756.

\bibitem[Shi]{onthewittvectorsofperfectringsinpositivecharacteristic}Kazuma \textsc{Shimomoto}, \textit{On the Witt Vectors of Perfect Rings in Positive Characteristic}. Communications in Algebra \textbf{43} n°12 (2015), pp. 5328--5342.

\bibitem[Sou]{anneauxuniform}Jean-Pierre \textsc{Soublin}, \textit{Anneaux uniformément cohérents}. Comptes rendus hebdoma-daires des séances de l'Académie des Sciences Série A \textbf{267} n°5 (1968), pp. 205--208.

\bibitem[Sta]{stacksproject}The \textsc{Stacks Project} Authors, \textit{Stacks Project}. \texttt{https://stacks.math.columbia.edu/} (2025).

\bibitem[Zin]{thedisplayofaformal}Thomas \textsc{Zink}, \textit{The display of a formal $p$-divisible group}. Astérisque \textbf{278} (2002), pp. 127--248.
\end{thebibliography}
\end{document}